\documentclass[12pt]{article}%
\usepackage[applemac]{inputenc}
\usepackage{amsmath,amssymb,fullpage}
\usepackage{amsfonts}
\usepackage{amsmath,amssymb}
\usepackage[applemac]{inputenc}
\usepackage{amsmath,amssymb,fullpage}
\usepackage{color}
\usepackage{color, colortbl}
\usepackage{amsmath}
\usepackage{amssymb}
\usepackage{graphicx}
\usepackage{tikz}
\usepackage{geometry}
\usetikzlibrary{arrows}
\usepackage{amsfonts}
\usepackage{amsmath,amssymb}
\usepackage[applemac]{inputenc}
\usepackage{amsmath,amssymb,fullpage}
\usepackage{color}
\usepackage{amsmath}
\usepackage{amssymb}
\usepackage{graphicx}
\usepackage{tikz}
\newtheorem{theorem}{Theorem}[section]
\newtheorem{lemma}[theorem]{Lemma}

\newtheorem{definition}[theorem]{Definition}

\newtheorem{example}[theorem]{Example}

\newtheorem{problem}[theorem]{Problem}
\newtheorem{remark}[theorem]{Remark}

\usetikzlibrary{arrows,shapes,automata,petri}

\newcommand{\dproof}{\noindent {Proof.} \quad}
\newcommand{\fproof}{\hfill $\square$ \bigskip}

\numberwithin{equation}{section}

\def\RB{\mathbb{R}}

\def\AC{\mathcal{A}}

\definecolor{LightCyan}{rgb}{0.88,1,1}

\def\RR{{\mathbb{ R}}}

\def\EE{{\mathbb{ E}}}

\def\1B{\text{1\!\!I}}

\def\m{\mu}
\def\l{\langle}
\def\r{\rangle}
\def\<{\langle}
\def\>{\rangle}

\begin{document}

\title{Stochastic Fokker-Planck PIDE for conditional McKean-Vlasov jump diffusions and applications to optimal control}
\author{Nacira Agram$^{1}$ \& Bernt \O ksendal$^{2}$}
\date{2 September 2022}
\maketitle

\footnotetext[1]{Department of Mathematics, KTH Royal Institute of Technology 100 44, Stockholm, Sweden. \newline
Email: nacira@kth.se. Work supported by the Swedish Research Council grant (2020-04697).}

\footnotetext[2]{%
Department of Mathematics, University of Oslo, Norway. 
Email: oksendal@math.uio.no.}

\begin{abstract}
The purpose of this paper is to study optimal control of conditional McKean-Vlasov  (mean-field) stochastic differential equations with jumps (conditional McKean-Vlasov jump diffusions, for short). 
To this end, we first prove a stochastic Fokker-Planck equation for the conditional law of the solution of such equations. 

Combining this equation with the original state equation, we obtain a Markovian system for the state and its conditional law. Furthermore, we apply this to formulate an Hamilton-Jacobi-Bellman (HJB) equation for the optimal control of conditional McKean-Vlasov jump diffusions.

Then we study the situation when the law is absolutely continuous with respect to Lebesgue measure. In that case the Fokker-Planck equation reduces to a stochastic partial differential equation (SPDE)  for the Radon-Nikodym derivative of the conditional law.

Finally we apply these results to solve explicitly the following problems:
\begin{itemize}
    \item Linear-quadratic optimal control of conditional stochastic McKean-Vlasov jump diffusions.
    \item Optimal consumption from a cash flow modelled as a conditional stochastic McKean-Vlasov differential equation with jumps.  
\end{itemize}
\end{abstract}



\textbf{Keywords :}  Jump diffusion; common noise; conditional McKean-Vlasov differential equation; stochastic Fokker-Planck equation; optimal control; HJB equation.

\section{Introduction}
A conditional McKean-Vlasov equation is a stochastic differential equation (SDE) where the coefficients depend on both the state of the solution and its probability law, conditioned on some background noise, called common noise. In the unconditional case this type of equation was first studied by H. McKean in  \cite{M}. 

The first study of optimal control of such systems was done by Andersson \& Djehiche \cite {AD}, who introduced a stochastic maximum principle approach and solved a mean-variance portfolio selection problem. It was later extended to jump diffusion by many authors including Hafayed \cite{H}, and even to regime-switching for mean-field systems with jumps by Bayraktar \& Chakraborty \cite{BC}.
An early discussion of a stochastic Fokker-Planck equation for McKean-Vlasov systems with common noise, arising as  a limit of an increasing system of interacting particles, was done by Kolokoltsov \& Troeva in an unpublished paper \cite{KT}.
Buckdahn et al \cite{BLPR} prove that the expected value of a function of the solution of a mean-field stochastic differential equation at the terminal time satisfies a non-local partial differential equation (PDE) of mean-field type. They are however, not applying this to optimal control. 
Bensoussan et al \cite {BHL} work directly on a deterministic PDE of Fokker-Planck type, which they assume is satisfied by the law of a solution of a corresponding mean-field SDE. Then they study optimal control and Nash equilibria of games for such deterministic systems, by means of an HJB equation. We refer also to Lauri\`ere and Pirroneau \cite {LP,LP1}, where models with constant volatility have been considered.
Guo et al \cite{GPW} established It\^o's formula along a flow of probability measures that enables derivation of dynamic programming equations and verification theorems for controlled McKean-Vlasov jump diffusions. 
Miller \& Pham \cite{MP} have studied linear-quadratic McKean-Vlasov stochastic differential games. However, the mean-field appears as a conditional expectation of the state and the control but not the conditional law as in our general setting. Moreover, to solve the optimal control problem, the authors used a weak submartingale optimality principle. For existence and uniqueness of the  solution of McKean-Vlasov SDEs, and the associated Fokker-Planck equation, we refer to Jourdain et al \cite{JMW}, and to Bogachev et al \cite {BKRS} and Barbu \& R\" ockner \cite {BR,BR1}. 

The paper which seems to be closest to our paper is Pham \& Wei \cite{PW}. They derive a dynamic programming principle for conditional mean-field systems, and use this to prove that the value function is a viscosity solution of an associated HJB equation.

Our paper differs from the above papers in several ways:
\begin{itemize}
    \item 
We include jumps in the system. As far as we know none of the related papers in the literature are dealing with conditional McKean-Vlasov  jump diffusions.
\item
Our methods are different. We use Fourier transform of measures to derive a general Fokker-Planck equation in the sense of distributions for the law of the solution of a conditional McKean-Vlasov equation. 
\item
As a result of this, we obtain a combined Markovian stochastic differential equation for the state and its law process. Then we use this to prove an HJB equation for optimal control of McKean-Vlasov jump diffusions. We derive \emph{sufficient} conditions of optimality (a \emph{verification} theorem).
\item
Our method allows us to obtain explicit solutions of some optimal control problems for McKean-Vlasov jump diffusions.
\end{itemize}
The paper is organised as follows: In Section 2 we recall some preliminaries that will be used throughout this work. In Section 3 we prove a Fokker-Planck equation for the conditional law of the solution of a MacKean-Vlasov jump diffusion and in Section 4 we study the optimal control of conditional McKean-Vlasov jump diffusions by means of a Hamilton-Jacobi-Bellman (HJB) equation. Section 5 presents the Fokker-Planck equation in the absolutely continuous case and gives an HJB equation in that situation.
Finally, in Section 6 we illustrate our theory by solving explicitly a linear-quadratic control problem for conditional McKean-Vlasov jump diffusions and a problem of optimal consumption from a cash flow modelled by a McKean-Vlasov jump diffusion with common noise.

\section{Preliminaries}
We now recall some basic concepts and background results:
\subsection{Radon measures}
 A Radon measure on $\mathbb{R}^d$ is a Borel measure which is finite on compact sets, outer regular on all Borel sets and inner regular on all open sets. In particular, all Borel probability measures on $\mathbb{R}^d$ are Radon measures.\\
In the following, we let $\mathbb{M}_0$ be the set of deterministic Radon measures and we let $C_0(\mathbb{R}^d)$ be the uniform closure of the space $C_c(\mathbb{R}^d)$ of continuous functions with compact support. If we equip $\mathbb{M}_0$ with the total variation norm $||\mu||:=|\mu|(\mathbb{R}^d)$, then $\mathbb{M}_0$ becomes a Banach space, and it is the dual of $C_0(\mathbb{R}^d)$. See  Chapter 7 in Folland \cite{F} for more information.\\
If $\mu \in \mathbb{M}_0$ is a finite measure, we define
\begin{equation}\label{mu^}
\widehat{\mu}(y):=F[\mu](y):=\int_{\mathbb{R}^d}e^{-ixy} \mu(dx); \quad y \in \mathbb{R}^d
\end{equation}
to be the Fourier transform of $\mu $ at $y$. \\
In particular, if $\mu(dx)$ is absolutely continuous with respect to Lebesgue measure $dx$ with Radon-Nikodym-derivative $m(x)=\frac{\mu(dx)}{dx}$, so that $\mu(dx)=m(x) dx$ with $m \in L^1(\RR^d)$, we define
the Fourier transform
of $m$ at $y$, denoted by $\widehat{m}(y)$ or $F[m](y)$, by
\begin{align*}
F[m](y)=\widehat{m}(y)=\int_{\mathbb{R}^d} e^{-ixy} m(x) dx; \quad y \in \mathbb{R}^d.
\end{align*}
Here and in the following, we define
$$
L^p(\RR^d)=\{f:\RR^d \mapsto \RR ; \int_{{\RR}^{d}} |f(x)|^pdx <  \infty  \}, \quad p \in (0,\infty),
$$
and if $\mathcal{X},\mathcal{Y},\mathcal{Z}$ are Banach spaces, we let $C^{1,2,2}(\mathcal{X}\times\mathcal{Y}\times \mathcal{Z})$ denote the set of real functions on $\mathcal{X}\times\mathcal{Y}\times \mathcal{Z}$ which are continuously differentiable with respect to the first variable, and  twice continuously differentiable with respect to the two last variables.





\subsection{Schwartz space of tempered distributions}
Let $\mathcal{S}=\mathcal{S}(\mathbb{R}^d)$\label{simb-028} be \textit{the
Schwartz space of rapidly decreasing smooth real
functions} 
on $\mathbb{R}^d$. The space $\mathcal{S}=\mathcal{S}(\mathbb{R}^d)$ is a
Fr\'echet space with respect to the family of seminorms:\label{simb-029} 
\begin{equation*}
\Vert f \Vert_{k,\alpha} := \sup_{x \in \mathbb{R}^d}\big\{ (1+|x|^k) \vert
\partial^\alpha f(x)\vert \big\},
\end{equation*}
where $k = 0,1,...$, $\alpha=(\alpha_1,...,\alpha_d)$ is a multi-index with $%
\alpha_j= 0,1,...$ $(j=1,...,d)$ and\label{simb-030} 
\begin{equation*}
\partial^\alpha f := \frac{\partial^{|\alpha|} f}{\partial
x_1^{\alpha_1}\cdots \partial x_d^{\alpha_d}}
\end{equation*}
for $|\alpha|=\alpha_1+ ... +\alpha_d$. \\
Let $\mathcal{S}^{\prime }=\mathcal{S}^{\prime }(\mathbb{R}^{d})$\label%
{simb-031} be its dual, called the space of \emph{tempered distributions}. 
\index{tempered distributions} 
If $\Phi \in \mathcal{S}^{\prime }$ and $f \in \mathcal{%
S}$ we let \label{simb-033} 
\begin{equation}
\Phi (f) \text{ or } \langle \Phi ,f \rangle  \label{3.1}
\end{equation}%
denote the action of $\Phi $ on $f $. For example, if $\Phi =m$ is a
measure on $\mathbb{R}^{d}$ then 
\begin{equation*}
\langle \Phi,f \rangle =\int\limits_{\mathbb{R}^{d}} f(x)dm(x),
\end{equation*}%
and, in particular, if this measure $m$ is concentrated on $x_{0}\in \mathbb{%
R}^{d}$, then 
\begin{equation*}
\langle \Phi , f \rangle = f(x_{0})
\end{equation*}
is evaluation of $f $ at $x_{0}\in \mathbb{R}^{d}$.\newline
Other examples include 
\begin{equation*}
\left\langle \Phi , f \right\rangle = f^{\prime }(x_{1}),
\end{equation*}%
i.e. $\Phi $ takes the derivative of $f$ at a point $x_{1}$. More generally, 
$\left\langle \Phi , f \right\rangle = f^{(k)}(x_{k})$
i.e. $\Phi $ takes the $k$'th derivative at the point $x_{k}$,
or linear combinations of the above.\newline
If $\Phi \in \mathcal{S}'$ we define its Fourier transform $\widehat{\Phi} \in \mathcal{S}'$
by the identity
\begin{align*}
\langle \widehat{\Phi}, f \rangle = \langle \Phi, \widehat{f} \rangle; \quad f \in \mathcal{S}.
\end{align*}
The partial derivative with respect to $x_k$ of a tempered distribution $\Phi$ is defined by 
\begin{align*}
&\langle \tfrac{\partial}{\partial x_k} \Phi, f \rangle = -\langle \Phi, \tfrac{\partial}{\partial x_k} f \rangle; \quad \phi \in \mathcal{S}^d.\\
&\text{ More generally, }\\
&\langle \partial^{\alpha} \Phi, f \rangle = (-1)^{|\alpha|} \langle \Phi, \partial^{\alpha} f \rangle; \quad \phi \in \mathcal{S}^d.
\end{align*}
We refer to Chapter 8 in Folland \cite{F} for more information.\\

\section {Stochastic Fokker-Planck equation for the conditional McKean-Vlasov jump diffusion}

Let $X(t)=X_t \in \RR^d$ be a mean-field stochastic differential equation with jumps, from now on called a \emph{McKean-Vlasov jump diffusion}, of the form (using matrix notation),
\begin{align}\label{MK}
dX(t)& =\alpha(t,X(t),\mu_t)dt+\beta(t,X(t),\mu_{t})dB(t)+\int_{
\mathbb{R}^d
}\gamma(t,X(t^-),\mu_{t^{-}},\zeta)\widetilde{N}(dt,d\zeta), \nonumber\\ 
X(0) &=x\in\mathbb{R}^d,
\end{align}
where $B(t) \in \RR^m = \RR^{m\times 1}, \widetilde{N} \in \RR^k=\RR^{k\times 1}$ are, respectively, an $m$-dimensional  Brownian motion and a $k$-dimensional compensated Poisson random measure on a filtered probability space $(\Omega,\mathcal{F},\mathbb{F}=\{\mathcal{F}_t\}_{t\geq 0},P)$. \\

For convenience, we assume that for all $\ell; 1 \leq \ell \leq k,$ the L\'evy measure of $N_{\ell}$, denoted by $\nu_{\ell}$, satisfies the condition
$\int_{\RR} \zeta^2 \nu_{\ell}(d\zeta) < \infty$, which means that $N_{\ell}$ does not have many big jumps (but $N_{\ell}$ may still have infinite total variation near $0$). This assumption allows us to use the version of the It\^o formula for jump diffusion given in Theorem 1.16 in \cite{OS}.\\
Here $\mu_t=\mu_t(\omega)=\mu_t(dx,\omega) $ is the \emph{conditional law} of $X(t)$ given the filtration $\mathcal{F}_t^{(1)}$ generated by the first component $B_1$ of the $m$-dimensional Brownian motion $B$. More precisely, we consider the following model:
\begin{definition}
We assume that $m \geq 2$ and we fix one of the Brownian motions, say $B_1=B_1(t,\omega)$, with filtration $\{\mathcal{F}_t^{(1)}\}_{t\geq 0}$. We define $\mu_t=\mu_t(\omega, dx)$ to be regular conditional distribution of $X(t)$ given $\mathcal{F}_t^{(1)}$. This means that $\mu_t(\omega,dx) $ is a Borel probability measure on $\RR^d$ for all $t \in [0,T],\omega \in \Omega$ and
\begin{equation} \label{cond}
\int_{\mathbb{R}^n} g(x) \mu_t(dx,\omega)= \EE[g(X(t)) | \mathcal{F}_t^{(1)}](\omega)
\end{equation}
 for all functions $g$ such that $\EE[ |g(X(t)) |] < \infty$.
\end{definition}

\begin{remark}
Heuristically, the equation \eqref{MK} models a mean-field system which is subject to what is called a "common noise" coming from the Brownian motion $B_1(t)$, which can be observed and is influencing the dynamics of the system, in the sense that the probability law $\mathcal{L}(X(t))$ of the state $X(t)$ is replaced by the conditional law $\mathcal{L}(X(t) | \mathcal{F}_t^{(1)})$. See e.g. Pham \& Wei \cite{PW}.
\end{remark}
Note that $\mu_t \in \mathbb{M}_0$ for all fixed $t, \omega$. From now on we let $\mathbb{M}$ denote all random measures $\lambda(dx,\omega)$ which are Radon measures for each $\omega$. \\ 
Throughout this paper, we assume the following assumptions on the coefficients:\\
$\alpha(t,x,\mu):[0,T]\times \mathbb{R}^d \times \mathbb{M}\rightarrow \mathbb{R}^d, \beta(t,x,\mu):[0,T]\times \mathbb{R}^d \times \mathbb{M}\rightarrow \mathbb{R}^{d \times m}$ and $\gamma(t,x,\mu,\zeta):[0,T]\times \mathbb{R}^d\times \mathbb{M}\times \mathbb{R}^d \rightarrow \mathbb{R}^{d \times k}$ are bounded deterministic functions for all $x,\mu, \zeta$, 
and that $\alpha, \beta, \gamma$ are continuous with respect to $t$ and $x$ for all $\mu,\zeta$. \newline

It was shown in Kurtz and Xiong \cite{KX} that under suitable assumptions on the coefficients there exists a unique solution of equation  \eqref{MK} with $\mathbb{E}[X^2(t)] < \infty$ for all $t$ in the continuous case (when $\nu=0$). Extension to jumps followed by the same arguments in \cite{KX}. \newline


In this section we will prove a stochastic Fokker-Planck equation for the conditional distribution of the McKean-Vlasov jump diffusion.\\

First, we recall some notations.\\
For fixed $t,\mu,\zeta$ and $\ell=1,2,... k$ we write for simplicity $\gamma^{(\ell)}=\gamma{(\ell)}(t,x,\mu,\zeta)$ for column number $\ell$ of the $d \times k$-matrix $\gamma$. For given $\mu \in \mathbb{M}$ the map
$$g \mapsto \int_{\mathbb{R}^d} g(x+\gamma^{(\ell)}) \mu(dx)$$
 is a bounded linear map on $C_0(\mathbb{R}^d)$. 
Therefore there is a unique measure $\mu^{(\gamma^{(\ell)})} \in \mathbb{M}$ such that 
\begin{align}\label{4.1}
\langle \mu^{(\gamma^{(\ell)})},g \rangle:=\int_{\mathbb{R}^d} g(x) \mu^{(\gamma^{(\ell)})}(dx)=\int_{\mathbb{R}^d}g(x+\gamma^{(\ell)}) \mu(dx), \text{ for all } g \in C_0(\mathbb{R}^d).
\end{align}
We call $\mu^{(\gamma^{(\ell)})}$ the $\gamma^{(\ell)}$-shift of $\mu$. \\ Note that $\mu^{(\gamma^{(\ell)})}$ is positive
 and absolutely continuous with respect to $\mu$.\\
 
In the following, we use $D_j, D_{n,j}$
to denote $\frac{\partial }{\partial x_j}$  and  $\frac{\partial^2}{\partial x_n \partial x_j}$ for notational simplicity, in the sense of distributions. The purpose of this section is to prove the following:

\begin{theorem}{(Conditional stochastic Fokker-Planck equation)} \label{FP1}\\
Let $X(t)$ be as in \eqref{MK} and let $\mu_t=\mu_t(dx,\omega)$ be the regular conditional distribution of $X(t)$ given $\mathcal{F}_t^{(1)}$. Then $\mu_t$ satisfies the following SPDE (in the sense of distributions):
 
\begin{align} \label{FPmu}
d\mu _{t} =A_0^{*} \mu_t dt + A_1^{*}\mu_t dB_1(t);   \quad \mu_0=\mathcal{L}(X(0)),
\end{align}
where 
$A_0^{*}$ is the integro-differential operator
\begin{align}
A_0^{*}\mu&= -\sum_{j=1}^d D_j [\alpha_j \mu] +\frac{1}{2}\sum_{n,j=1}^d D_{n,j}[(\beta \beta^{(T)})_{n,j} \mu] \nonumber\\
&+\sum_{\ell=1}^k  \int_{\mathbb{R}}\Big\{\mu^{(\gamma^{(\ell)})}-\mu+\sum_{j=1}^d D_j[\gamma_j^{(\ell)}(s,\cdot,\zeta)\mu]\Big\} \nu_{\ell} \left( d\zeta \right) \label{A0*}
\end{align}
and
\begin{align}
A_1^{*}\mu= - \sum_{j=1}^d D_j[\beta_{1,j} \mu], \label{A1*}
\end{align}
were $\beta^{(T)}$ denotes the transposed of the $d \times m$ - matrix $\beta=\big[\beta_{j,k}\big]_{1\leq j \leq d,1 \leq k \leq m}$  and $\gamma^{(\ell)}$ is column numer $\ell$ of the matrix $\gamma$.

\end{theorem}

\dproof 
Choose $\psi \in C^{2}\left( \mathbb{R}^d\right) $ with bounded derivatives, and with values in the complex plane $\mathbb{C}$. Then since $B_1$ is independent of the the other Brownian motions (and of the random measures $N_{\ell}$ ), we get by the It\^{o}
formula for jump diffusions (see e.g. Theorem 1.16 in \cite{OS}):
\begin{align}
&\mathbb{E}\Big[ \psi \left( X_{t}\right)| \mathcal{F}_t^{(1)}\Big]-\psi(x)\nonumber \\
&= \int_{0}^{t}\EE[A_0\psi \left( X_{s}\right) | \mathcal{F}_t^{(1)}]ds
+\int_0^{t}\EE[ A_1 \psi(X_s) |\mathcal{F}_t^{(1)}] dB_1(s)\nonumber\\
&= \EE\Big[\Big(\int_{0}^{t}\EE[A_0\psi \left( X_{s}\right) | \mathcal{F}_s^{(1)}]ds
+\int_0^{t}\EE[ A_1 \psi(X_s) |\mathcal{F}_s^{(1)}] dB_1(s)\Big) |\mathcal{F}_t^{(1)}\Big]\nonumber\\
&= \int_{0}^{t}\EE[A_0\psi \left( X_{s}\right) | \mathcal{F}_s^{(1)}]ds
+\int_0^{t}\EE[ A_1\psi(X_s) |\mathcal{F}_s^{(1)}] dB_1(s), \label{3.6}
\end{align}
where
\small
\begin{align*}
A_0\psi \left( X_s \right) &=\sum_{j=1}^d \alpha_j \left( s,X_s,\mu_s\right) \frac{\partial \psi}{\partial x_j}\left( X_{s}\right)+\frac{1}{2}\sum_{n,j=1}^d (\beta \beta^{T})_{n,j} (s,X_s,\mu_s)\frac{\partial ^2 \psi}{\partial x_n \partial x_j} \left( X_{s}\right) \\
&+\sum_{\ell =1}^k \int_{\mathbb{R}}\Big\{ \psi \left( X_{s}+\gamma^{(\ell)} \left( s,X_s,\mu_s,\zeta
\right) \right) -\psi \left( X_{s}\right)\\& -\sum_{j=1}^d \frac{\partial \psi}{\partial x_j} \left(
X_{s}\right) \gamma_j^{(\ell)} \left( s,X_s,\mu_s,\zeta \right) \Big\} \nu_{\ell} \left( d\zeta
\right) ,
\end{align*}
and 
\begin{align*}
A_1\psi(X_s)= \sum_{j=1}^d\frac{\partial \psi}{\partial x_j} (X_s)\beta_{1,j} (s,X_s,\mu_s),
\end{align*}
where  $\nu_{\ell} \left( \cdot \right) $ is the L\'{e}vy measure of $N_{\ell}\left( \cdot,\cdot\right). $\\
In particular, choosing, with $i=\sqrt{-1}$,
\[
\psi \left(x \right) =\psi _{y}\left( x\right) =e^{-iyx};\quad y,x\in 
\mathbb{R}^d,
\]
we get
\small
\begin{align*} 
&A_0\psi _{y}\left( X_{s}\right) \nonumber\\
& =\left( -i\sum_{i=1}^d y_i \alpha_i \left( s,X_{s},\mu_s\right) -%
\frac{1}{2}\sum_{n,j=1}^d y_n y_j (\beta \beta^{(T)})_{n.j} \left( s,X_{s},\mu_s\right) \right.  \nonumber\\
&\left. +\sum_{\ell=1}^k \int_{\mathbb{R}}\left\{ \exp \left( -iy\gamma^{(\ell)} \left(
s,X_{s},\mu_s,\zeta \right) \right) -1+i \sum_{j=1}^d y_i \gamma_j^{(\ell)} \left( s,X_{s},\mu_s,\zeta \right)
\right\} \nu_{\ell} \left( d\zeta \right) \right) e^{-iyX_{s}}.\label{3.4}
\end{align*}
and
\begin{align}
A_1 \psi(X_s)= -i \sum_{j=1}^d y_j\beta_{1,j}(s,X_s,\mu_s). 
\end{align}
In general we have (see \eqref{cond})
\[
\mathbb{E}\left[ g\left( X_{s}\right) e^{-iyX_{s}} | \mathcal{F}_s^{(1)}\right] =\int_{\mathbb{R}^d
}g\left( x\right) e^{-iyx}\mu _{s}\left( dx\right) =F\left[ g\left( \cdot\right)
\mu _{s}(\cdot)\right] \left( y\right).
\]
Therefore, we get
\begin{align}
\EE[ e^{-i\gamma(s,X_s,\zeta)}e^{-iyX_s} | \mathcal{F}_s^{(1)}] &= \int_{\RR^d} e^{-iy\gamma(s,x,\zeta)} e^{-ixy}\mu_s(dx)
=\int_{\RR^d} e^{-iy(x+\gamma(s,x,\zeta))} \mu_s(dx)\nonumber\\
&=\int_{\RR^d} e^{-iyx} \mu^{(\gamma)}(dx)=F[\mu^{(\gamma)}(\cdot)](y),\label{3.5}
\end{align}
where 
$\mu_s^{(\gamma)}(\cdot)$ is the $\gamma$-shift of $\mu_s$. Recall that if $w \in \mathcal{S}'$, using the notation $\frac{\partial}{dx_j}w\left( t,x\right) =:D_jw\left( t,x\right) ,$ and similarly with higher order derivatives, we have, in the sense of distributions,
\[
F\left[ D_jw\left( t,\cdot  \right) \right] ( y)
=iy_j F\left[ w\left( t,\cdot\right) \right] \left( y\right). 
\]
Therefore, 
\begin{align}\label{D-alpha}
iy_jF[\alpha(s. \cdot)\mu_s](y)&= F[D_j(\alpha(s,\cdot)\mu_s)](y)\\\label{D-beta}
-\ y_ny_jF[\beta \beta^{T}(s,\cdot)\mu_s](y)&=F[D_{n,j}(\beta \beta^{T}(s,\cdot) \mu_s)](y).
\end{align}
Applying this and \eqref{3.5} to \eqref{3.4}, we get 
\begin{align*}
&\EE[A_0\psi_y \left( X_{s}\right) | \mathcal{F}_s^{(1)}]= \int_{\RR^d} \Big(-i\sum_{j=1}^d y_j \alpha_j \left( s,x,\mu_s\right) -\frac{1}{2}\sum_{n,j=1}^d y_n y_j (\beta \beta^{(T)})_{n,j} \left( s,x,\mu_s\right)   \nonumber\\
&+\sum_{\ell=1}^k \int_{\mathbb{R}}\Big\{ \exp \left( -iy\gamma^{(\ell)} \left(
s,x,\mu_s,\zeta \right) \right) -1+i \sum_{j=1}^d y_j\gamma_j^{(\ell)} \left( s,x,\mu_s,\zeta \right)
\} \nu_{\ell} \left( d\zeta \right)\Big)e^{-iyx}\mu_s(dx) \nonumber\\
&=-i\sum_{j=1}^d y_j F[\alpha_j \mu_s](y) -\frac{1}{2}\sum_{n,j=1}^d y_n y_j F[(\beta \beta^{(T)})_{n.j} \mu_s]  \nonumber\\
&+\sum_{\ell=1}^k F\Big[ \int_{\RR}\Big\{ \exp \left( -iy\gamma^{(\ell)} \left(
s,x,\mu_s,\zeta \right) \right) -1+i \sum_{j=1}^d y_j \gamma_j^{(\ell)} \left( s,x,\mu_s,\zeta \right)
\Big\} \nu_{\ell} \left( d\zeta \right) \mu_s\Big](y)\nonumber\\
&=F\Big[ -\sum_{j=1}^d D_j [\alpha_j \mu_s] +\frac{1}{2}\sum_{n,j=1}^d D_{n,j}[(\beta \beta^{(T)})_{n.j} \mu_s] \nonumber\\
&+\sum_{\ell=1}^k  \int_{\mathbb{R}}\Big\{\mu_s^{(\gamma^{(\ell)})}-\mu_s+\sum_{j=1}^d D_j[\gamma_j^{(\ell)}(s,\cdot,\zeta)\mu_s]
\Big\} \nu_{\ell} \left( d\zeta \right)\Big](y)\nonumber\\
&=F[A_0^{*} \mu_s](y),
\end{align*}
where $A_0^{*}$ is the integro-differential operator 
\begin{align*}
A_0^{*}\mu&= -\sum_{j=1}^d D_j [\alpha_j \mu] +\frac{1}{2}\sum_{n,j=1}^d D_{n,j}[(\beta \beta^{(T)})_{n,j} \mu] \nonumber\\
&+\sum_{\ell=1}^k  \int_{\mathbb{R}}\{\mu^{(\gamma^{(\ell)})}-\mu+\sum_{j=1}^d D_j[\gamma_j^{(\ell)}(s,\cdot,\zeta)\mu]\} \nu_{\ell} \left( d\zeta \right)
\end{align*}
Note that $A^{*}\mu_s$ exists in $\mathcal{S}^{\prime}$.\\
Similarly, we get
\begin{align*}
\EE[ A_1 \psi(X_s) |\mathcal{F}_s^{(1)}]&= \int_{\RR} -i \sum_{j=1}^d y_j\beta_{1,j}(s,x,\mu_s)e^{-iyx} \mu_s(dx)\nonumber\\
&=F[-i \sum_{j=1}^d y_j\beta_{1,j}(s,x,\mu_s)\mu_s] =F[-\sum_{j=1}^d D_j[\beta_{1,j} \mu_s]](y) \nonumber\\
&=F[ A_1^{*}\mu_s](y),
\end{align*}
where $A_1^{*}$ is the operator
\begin{align*}
A_1^{*}\mu_s= - \sum_{j=1}^d D_j[\beta_{1,j} \mu_s].
\end{align*}
Hence
\begin{align}
&\mathbb{E}\Big[ \psi \left( X_{t}\right) | \mathcal{F}_t^{(1)}\Big]\nonumber \\
&=\psi(x)+  \int_{0}^{t}\EE[A_0\psi \left( X_{s}\right) | \mathcal{F}_t^{(1)}]ds
+\int_0^{t}\EE[ A_1(s) |\mathcal{F}_t^{(1)}] dB_1(s)\nonumber\\
&=\psi(x)+ \EE[ \int_{0}^{t}\EE[A_0\psi \left( X_{s}\right) | \mathcal{F}_s^{(1)}]ds
+\int_0^{t}\EE[ A_1(s) |\mathcal{F}_s^{(1)}] dB_1(s)|\mathcal{F}_t^{(1)}]\nonumber\\
&= \psi(x)+\EE[\int_0^{t}F [A_0^{*}\mu_s](y)ds +\int_0^{t} F[A_1^{*}\mu_s](y) dB_1(s)|\mathcal{F}_t^{(1)}]\nonumber\\
&= \psi(x)+\int_0^{t}F [A_0^{*}\mu_s](y)ds +\int_0^{t} F[A_1^{*}\mu_s](y) dB_1(s).\label{3.15}
\end{align}
On the other hand,
\begin{align}
&\mathbb{E}\Big[ \psi \left( X_{t}\right)  | \mathcal{F}_t^{(1)}\Big]\nonumber =\EE[e^{-iyX_{t} }- e^{-iyX_0} |\mathcal{F}_t^{(1)}]\nonumber\\
&= \EE[(\EE[e^{-iyX_{t} }|\mathcal{F}_{t}^{(1)}]- e^{-iyX_0}) |\mathcal{F}_t^{(1)}]=\EE[\big(\widehat{\mu}_{t}(y)-\widehat{\mu}_0 (y)\big) |\mathcal{F}_t^{(1)}]=\widehat {\mu}_t(y) -\widehat{\mu}_0(y).\label{3.16}
\end{align}
Combining \eqref{3.15} and \eqref{3.16}, we get
\begin{align*}
\widehat{\mu}_{t}(y)-\widehat{\mu}_0(y)=\int_0^{t}F [A_0^{*}\mu_s](y)ds +\int_0^{t} F[A_1^{*}\mu_s](y) dB_1(s).
\end{align*}
Since the Fourier transform of a distribution determines the distribution uniquely, we deduce that
\begin{align*}
\mu_{t}-\mu_0=\int_0^{t} A_0^{*}\mu_s ds +\int_0^{t}  A_1^{*}\mu_s dB_1(s),
\end{align*}
or, in differential form,
\begin{align*}
d\mu_t=A_0^{*}\mu_t dt + A_1^{*}\mu_t dB_1(t); \quad \mu_0=\mathcal{L}(X(0)),
\end{align*}
as claimed.
\fproof 

\begin{remark}
\begin{itemize}
    \item In Theorem 3.3 we only prove existence of solution of equation \eqref{FPmu}. The uniqueness of solution of equation \eqref{FPmu} is studied by \cite{CG} in the case without jumps. They prove that under certain smoothness and Lipschitz conditions on the coefficients, the solution is unique, see Theorem 5.4 in  \cite{CG}. By inspecting the proof, we see that a similar result can be obtained in the jump case, provided corresponding assumptions are made on the jump coefficient. In the present paper, we will assume that the solution of \eqref{FPmu} is unique.
    
    \item Note that in the unconditional case, where the conditional law $\mathcal{L}(X(t) |\mathcal{F}_t^{(1)})$ is replaced by the law $\mathcal{L}(X(t))$, is a special case, which can be obtained from the above by putting $\beta_1=0$. In that case we get $A_1^{*}=0$, and the Fokker-Planck equation becomes a deterministic PDE. Similar remarks apply to the following sections, including the HJB equation.
   
\end{itemize}

\end{remark}

\begin{section} {An HJB equation for optimal control of conditional McKean-Vlasov jump diffusions (I)}
We now apply the results obtained in Section 3 to derive an HJB equation for optimal control of McKean-Vlasov jump diffusions.

Consider a controlled version of the system \eqref{MK}, in which we have introduced a control process $u=\{u(t),t\in[0,T] \}$, i.e. 

\begin{align} \label{MV3}
dX(t)&=dX^{(u)}(t)=\alpha(t,X(t),\mu_t,u(t))dt + \beta(t,X(t),\mu_t,u(t)) dB(t)\\
&+\int_{\mathbb{R}^k} \gamma(t,X(t^-), \mu_t,u(t), \zeta) \widetilde{N}(dt,d\zeta);\quad t>0,\nonumber\\
X(0)&= x \in \mathbb{R}^d. \nonumber
\end{align}
As before $\mu_t=\mu_t(dx,\omega)=\mathcal{L}(X(t)|\mathcal{F}_t^{(1)})$ denotes the conditional distribution of $X(t)$ given $\mathcal{F}_t^{(1)}$.\\
Let $\mathbb{U}$ denote the set of possible control values.
We assume that the coefficients\\ $\alpha(t,x,\mu,u):[0,T]\times \mathbb{R}^d \times \mathbb{M}\times \mathbb{U}\rightarrow \mathbb{R}^d, \beta(t,x,\mu):[0,T]\times \mathbb{R}^d \times \mathbb{M}\times \mathbb{U}\rightarrow \mathbb{R}^{d \times m}$ and $\gamma(t,x,\mu,\zeta,u):[0,T]\times \mathbb{R}^d\times \mathbb{M}\times \mathbb{R} ^k\times \mathbb{U}\rightarrow \mathbb{R}^{d \times k}$ are bounded deterministic functions for all $x,\mu, \zeta,u$, 
and that $\alpha, \beta, \gamma$ are continuous with respect to $t$ and $x$ for all $\mu,\zeta$. \\
We say that $u$ is \emph{admissible} if $u$ is Markovian, 
i.e. $u$ has the form $u(t)=u_0(t,X(t),\mu_t)$ for some function $u_0: [0,T]\times \RR^d \times \mathbb{M} \mapsto \RR$. \\
For all $u\in \mathbb{U}$, existence and unique of a solution $X^{(u)}$ of equation \eqref{MV3} such that 
\begin{align}
\EE[|X^{(u)}(t)|^2] < \infty  \text{ for all } t,
\end{align}
can be obtained in a same way as in Kurtz and Xiong \cite{KX}.

The set of admissible controls is denoted by $\mathcal{A}$. As is customary, for simplicity (and a slight abuse) of notation we do not distinguish in notation between $u_0$ and $u$ in the following.\\
From Theorem \ref{FP1}, we have the following stochastic Fokker-Planck equation for $\mu$:
\begin{align}\label{fp2}
d\mu_t &= A^{*,u}_{0} \mu_t dt + A_1^{*,u} \mu_t dB_1(t); \quad t>0,\\
\mu_0&=\mathcal{L}(X(0)),\nonumber
\end{align}
where the integro-differential operators $A_0^{*,u}, A_1^{*,u}$ are associated to the controlled process $u$, as follows (see \eqref{A0*},\eqref{A1*}):
\begin{align*}
A_0^{*,u}\mu&= A_0^{*}\mu=-\sum_{j=1}^d D_j [\alpha_j \mu] +\frac{1}{2}\sum_{n,j=1}^d D_{n,j}[(\beta \beta^{(T)})_{n,j} \mu] \nonumber\\
&+\sum_{\ell=1}^k  \int_{\mathbb{R}^k}\{\mu^{(\gamma^{(\ell)})}-\mu+\sum_{j=1}^d D_j[\gamma_j^{(\ell)}(s,\cdot,\zeta)\mu]\} \nu_{\ell} \left( d\zeta \right)
\end{align*}
and
\begin{align}
A_1^{*,u}\mu=A_1^{*}\mu= - \sum_{j=1}^d D_j[\beta_{1,j} \mu].
\end{align}
Combining \eqref{MV3} and \eqref{fp2} we can write the dynamics of the $d+2$-dimensional $[0,T] \times \mathbb{R}^d\times  \mathbb{M}$ - valued process $Y(t)=(Y_0(t),Y_1(t),Y_2(t))=(s+t, X(t),\mu_t)$ as follows:
\begin{align}\label{Y2}
dY(t)&=\left[\begin{array}{clcr}
dY_0(t)\\dY_1(t)\\dY_2(t)
\end{array} \right] 
=\left[\begin{array}{clcr}
dt \\dX(t)\\d\mu_t
\end{array} \right] \nonumber\\
&=\left[ \begin{array}{c}
1 \\ \alpha(Y(t),u(t)) \\A_0^{*,u}\mu_t
\end{array} \right] dt
+\left[\begin{array}{rc}
0_{1\times m} \\ \beta(Y(t),u(t)) \\ A_1^{*,u}\mu_t,0,0 ...,0
\end{array} \right]dB(t)\nonumber\\
&+ \int_{\mathbb{R}^d} \left[ \begin{array}{rc}
0_{1\times k}\\ \gamma(Y(t^-),u(t),\zeta) \\0_{1\times k}
\end{array}\right] \widetilde{N}(dt,d\zeta),
\end{align}
where we have used the shorthand notation
\begin{align*}
\alpha(Y(t),u(t))&=\alpha(Y_0(t),Y_1(t),Y_2(t),u(t))\\
\beta(Y(t),u(t))&= \beta(Y_0(t),Y_1(t),Y_2(t),u(t))\\
\gamma(Y(t^-),u(t),\zeta)&=\gamma(Y_0(t^-),Y_1(t^-),Y_2(t^-),u(t),\zeta).
\end{align*}
Let $f:\mathbb{R} \times \mathbb{R}^d \times \mathbb{M} \times \mathbb{U} \mapsto \mathbb{R}$ and $g: \mathbb{R}^d \times \mathbb{M} \mapsto \mathbb{R}$ be given functions such that
\begin{align*}
\mathbb{E}^{y}\Big[ \int_0^T |f(s+t, X(t),\mu_t,u(t))|dt  + |g(X(T),\mu_T)| \Big] < \infty \text{ for all } u\in \mathcal{A}.
\end{align*}
We introduce the performance functional:
\begin{align*}
&J_u(s, x,\mu)=J(y)\nonumber\\
&=\mathbb{E}^{y}\Big[ \int_0^T f(s+t, X(t),\mu_t,u(t))dt  + g(X(T),\mu_T) \Big], \quad u \in \mathcal{A},
\end{align*}
where $y=(y_0, y_1,y_2)$. Thus time starts at $y_0=s$, $X(t)$ starts at $x=y_1$ and $\mu_t$ starts at $\mu=y_2$.
Then we define the value function:
\begin{align*}
\Phi(y)=\Phi(s,x,\mu)= \sup_{u \in \mathcal{A}} J_u (s,x,\mu).
\end{align*}
In the following, if $\varphi=\varphi(s,x,\mu) \in C^{1,2,2}([0,T] \times \mathbb{R}^d \times \mathbb{M})$, then  $D_{\mu} \varphi=\nabla_{\mu} \varphi \in L(\mathbb{M},\RR)$ (the set of bounded linear functionals on $\mathbb{M}$) denotes the 
Fr\' echet derivative (gradient) of $\varphi$ with respect to $\mu \in \mathbb{M}$. Similarly $D_{\mu}^2 \varphi$ denotes the double derivative of $\varphi$ with respect to $\mu$. It is an element of $L(\mathbb{M},L(\mathbb{M},\RR))$. By the Riesz representation theorem it may be regarded as an element of $L(\mathbb{M} \times \mathbb{M},\RR)$ (the bounded linear functionals on $\mathbb{M} \times \mathbb{M}$). See Appendix for details.
\begin{remark}
\begin{itemize}
    \item  As we have explained in Section 2.1, by equipping the Radon space of measures with their total variation norm, it becomes a Banach space and we can use Fr\' echet derivative to represent differentiation with respect to a measure. However, if we consider measures on a Wasserstein metric space, then we restrict ourselves only to measures with finite second moments. Therefore we cannot easily differentiate with respect to a measure in this space. To do this, we would need to identify the space with a Hilbert space; to get what is called the Lions derivative. We think it is easier in the current paper to work with the Banach space $\mathbb{M}$.
    \item The control processes we consider are of feedback/Markovian (closed-loop) form with respect to both the state $X_t$ and its
conditional law $\mu_t$. The extension to controls of open-loop form would require a different approach e.g. a maximum principle approach, and is left for future work.
\end{itemize}

\end{remark}
We now give a sufficient condition (a verification theorem) for optimal control of such a problem. More precisely, we formulate an HJB equation such that if a smooth function $\varphi$ satisfies the HJB equation, it coincides with the value function $\Phi$:

\begin{theorem}\textbf{(An HJB equation for optimal control of conditional McKean-Vlasov jump diffusions (I))} \\
For $v \in \mathbb{U}$ define $G_v$ to be the following integro-differential operator, which is \emph{the generator of the process $Y(t)=(s+t,X(t),\mu_t) \in [0,T] \times \RR^d \times \mathbb{M}$},  given the control value $v$, defined as follows:
\begin{align*}
&G_v \varphi(s,x,\mu)= \frac{\partial \varphi}{\partial s} +\sum_{j=1}^d \alpha_j (s,x,\mu,v)\frac{\partial \varphi}{\partial x_j} + \langle \nabla_{\mu} \varphi, A_0^{*,v} \mu \rangle \nonumber\\
&+ \tfrac{1}{2}\sum_{j,n=1}^{d}  (\beta \beta^{T})_{j,n}(s,x,\mu,v) \frac{\partial ^2 \varphi}{\partial x_j \partial x_n} + \tfrac{1}{2}\sum_{j=1}^d \beta_{j,1}\frac{\partial}{\partial x_j}\langle\nabla_{\mu} \varphi,A_1^{*,v}\mu\rangle \nonumber\\
&+\tfrac{1}{2} 
\langle A_1^{*,v}\mu, \langle D_{\mu}^2 \varphi,A_1^{*,v}\mu\rangle \rangle \nonumber\\
&+\sum_{\ell =1}^k \int_{\mathbb{R}} \{ \varphi(s, x+\gamma^{(\ell)}, \mu)) - \varphi(s,x,\mu) 
-\sum_{j=1}^d\gamma_j^{(\ell)}  \tfrac{\partial}{\partial x_j} \varphi(s,x,\mu) \}\nu_{\ell}(d\zeta) ;\nonumber\\
& \varphi  \in C^{1,2,2}([0,T] \times \mathbb{R}^d \times \mathbb{M}).
\end{align*}
Here  $\gamma=\gamma(s,x,\mu,u,\zeta)$ and $\gamma^{(\ell)}$ is column number $\ell$ of the $d \times k$- matrix $\gamma$.
Suppose there exists a function $\widehat \varphi(s,x,\mu)\in C^{1,2,2}([0,T] \times \mathbb{R}^d \times \mathbb{M})$ and a Markov control $\widehat{u}=\widehat{u}(y) \in \AC$ such that, for all  $y$,

\begin{align}
\sup_{v \in \mathbb{U}} \Big\{f(y,v)+G_v\widehat \varphi(y) \Big \} &= f(y,\widehat{u}(y))+G_{\widehat{u}(y)}\widehat \varphi(y) = 0 \label{HJB1}\\
& \text{ and}\nonumber\\
&\widehat \varphi(T,x,\mu)=g(x,\mu)\nonumber.
\end{align}
Then $\widehat{u}$ is an optimal control and $\widehat \varphi=\Phi$.
\end{theorem}

\dproof
Choose  $\varphi \in C^{1,2,2}([0,T] \times \mathbb{R}^d \times \mathbb{M}) $ and let $u=u(y)\in \mathcal{A}$ be a Markov control for the Markov process $Y(t)=(s+t,X(t),\mu_t)$.  Then by the Dynkin formula, we have
\begin{align}
\EE^y[ \varphi(Y(T))] = \varphi(s,x,\mu) +\EE^y \Big[\int_0^T G_u\varphi(Y(t))dt\Big]. \label{Dynkin}
\end{align}
Suppose $\varphi$ satisfies the conditions
\begin{align}
f(y,v) +G_v\varphi(y) &\leq 0 \text{ for all } y=(s,x,\mu) \text{ and all } v, \label{hjb1} \\
&\text{ and }\nonumber\\
\varphi(T,x,\mu) &\geq g(x,\mu) \text{ for all } x,\mu.\label{hjb2}
\end{align}
Then by \eqref{Dynkin} we get
\begin{align*}
\EE^y[\varphi(Y(T))] \leq \varphi(s,x,\mu) -\EE^y \Big[\int_0^T f(Y(t),u(t))dt\Big]
\end{align*}
or
\begin{align}
\varphi(s,x,\mu) &\geq \EE^y\Big[ \int_0^T f(s+t,X(t), \mu_t,u(t)) dt + \varphi(T,X(T),\mu_T)\Big]\nonumber\\
&\geq \EE^y\Big[ \int_0^T f(s+t,X(t),\mu_t,u(t)) dt + g(X(T),\mu_T) \Big] \nonumber\\
&= J_u(s,x,\mu).\label{6.16}
\end{align}
Since this holds for all $u \in \mathcal{A}$, we deduce that
\begin{align}
\varphi(s,x,\mu) \geq \sup_{u \in \mathcal{A}} J_u(s,x,\mu)=\Phi(s,x,\mu).\label{6.17}
\end{align}
Now assume that $\varphi:=\widehat{\varphi}$ and $u:=\widehat{u} \in \mathcal{A}$ satisfy \eqref{hjb1} and \eqref{hjb2}. Then \eqref{6.16} holds with \emph{equality}, i.e.
\begin{align*}
\widehat{\varphi}(s,x,\mu)=J_{\widehat{u}} (s,x,\mu).
\end{align*}
We therefore obtain the following string of inequalities
\begin{align*}
\Phi(s,x,\mu) \leq \widehat{\varphi}(s,x,\mu) = J_{\widehat{u}} (s,x,\mu) \leq \sup_{u \in \mathcal{A}} J_{u}(s,x,\mu) = \Phi(s,x,\mu).
\end{align*}
Since the first term and the last term are the same, we have \emph{equality} everywhere in this string. Hence $\widehat{\varphi}=\Phi$ and $\widehat{u}$ is optimal.
\fproof
\begin{remark}
In the general case, when there is no smooth solution to the HJB equation  \eqref{HJB1}, then
its solution should be interpreted in the viscosity sense as in \cite{BIRS}, \cite{CGKPR}. However, to the best of our knowledge, viscosity solutions for McKean-Vlasov jump diffusions with common noise have not been studied yet. It is a topic for future research. 
\end{remark}

\end{section}

\section{The Fokker-Planck equation in the absolutely continuous case}
In this section we assume that $\mu _{t}\left( dx\right) <<dx$ for all $t>0$, and we put 
\begin{align}m\left( t,x\right) =\frac{\mu _{t}\left( dx\right) }{dx}, \text{ so that  } 
\mu _{t}\left( dx\right)=m\left( t,x\right) dx; \quad t > 0. \label{abs}
\end{align} 
We make the following observation:
\begin{lemma}
Assume that 
\begin{align}\label{gamma}
\gamma=\gamma(s,x,m,\zeta)=\gamma(s,m,\zeta) \text{  does not depend on } x. 
\end{align} 
Then
\begin{align}\label{5.3}
m^{(\gamma^{(\ell)})}(t,x) &=m(t,x-\gamma^{(\ell)});\quad \text{ for all } t,x.
\end{align}
\end{lemma}

\dproof
By \eqref{4.1} and a change of variable,
\begin{align}\label{mgamma}
\int_{\mathbb{R}^d} g(x) m^{(\gamma^{(\ell)})}(t,x) dx=\int_{\mathbb{R}^d}g(x+\gamma^{(\ell)}) m(t,x)dx =\int_{\mathbb{R}^d}g(x) m(t,x-\gamma^{(\ell)})dx 
\end{align}
for all $g \in C_0(\mathbb{R}^d)$.
\fproof\\

\begin{definition}
Let $K(dx)$ be the measure on $\RR^d$ defined by
\begin{align} \label{K} 
K(dx)=(1+|x|^2)dx,
\end{align}
and define 
$$L_K^1(\RR^d)=\{ f:\RR^d \mapsto \RR; ||f||_K:=\int_{\RR^d} |f(x)| K(dx) < \infty\}.$$
\end{definition}
Note that since $\EE[|X(t)|^2] < \infty$, we have $m(t,\cdot) \in L_K^1(\mathbb{R}^d)$ for all $t>0$. 

\begin{remark}
In this absolutely continuous setting, and if we consider probability measures only, the space $L_K^1(\RR^d)$ can be related to the Wasserstein metric space
 $\mathcal{P}_2(\RR^d)$. See e.g. Cardaliaguet \cite{C} and Lions \cite{LL}. See also Agram \cite{A}, who applies this to an optimal control problem for McKean-Vlasov equations with anticipating law.
\end{remark}
With a slight change of notation, we can write the McKean-Vlasov jump equation \eqref{MK} as a stochastic differential equation in the $\mathbb{R}^d$-valued process $X(t)=X(t,\omega); (t,\omega)\in \ [0,T] \times \Omega$ as follows:
\begin{align} \label{MV2}
dX(t)&= \alpha(t,X(t),m(t,\cdot))dt + \beta(t,X(t),m(t,\cdot)) dB(t)\\
&+\int_{\mathbb{R}^k} \gamma(t,X(t^-), m(t,\cdot),\zeta) \widetilde{N}(dt,d\zeta);\quad t>0,\nonumber\\
X(0)&= x \in \mathbb{R}^d,\nonumber
\end{align}
where the coefficients $\alpha(t,x,m):[0,T]\times \mathbb{R}^d \times L_K^1(\mathbb{R}^d)\rightarrow \mathbb{R}^d,\beta(t,x,m):[0,T]\times \mathbb{R}^d \times L_K^1(\mathbb{R}^d)\rightarrow \mathbb{R}^{d\times m}$ and $\gamma(t,x,m,\zeta):[0,T]\times \mathbb{R}^d \times L_K^1(\mathbb{R}^d)\times \mathbb{R}^k \rightarrow \mathbb{R}^{d\times k}$ are bounded. 
In the previous section we proved that
\begin{align*}
d\mu_t = A_0^{*}\mu_t dt + A_1^{*}\mu_t dB_1(t) \text{ in } \mathcal{S}'.
\end{align*}
If $\mu_t(dx)=m(t,x)dx$, this equation becomes the SPDE
\begin{align}
dm(t,x)= A_0^{*}m(t,x) dt + A_1^{*}m(t,x) dB_1(t), 
\end{align}
where the operators $A_0^{*}, A_1^{*}$ get the form (see \eqref{A0*},\eqref{A1*}):
\begin{align}
A_0^{*}m&= - \sum_{j=1}^d \frac{\partial}{\partial x_j}[\alpha_j m] +\frac{1}{2}\sum_{n,j=1}^d \frac{\partial^2}{\partial x_j \partial x_n}[(\beta \beta^{(T)})_{n,j} m] \nonumber\\
&+\sum_{\ell=1}^k  \int_{\mathbb{R}}\Big\{m^{(\gamma^{(\ell)})}-m+\sum_{j=1}^d \frac{\partial}{\partial x_j} [\gamma_j^{(\ell)}(s,\cdot,\zeta)m]\Big\} \nu_{\ell} \left( d\zeta \right), \label{A_0m}
\end{align}
and
\begin{align}
A_1^{*}m= - \sum_{j=1}^d \tfrac{\partial}{\partial x_j}[\beta_{1,j} m].\label{A_1m}
\end{align}
Here $m^{(\gamma^{(\ell)})}$ is the Radon Nikodym derivative of $\mu^{(\gamma^{(\ell)})}$ with respect to Lebesgue measure, i.e.,
\begin{align}\label{5.2}
m^{(\gamma^{(\ell)})}(t,x)= \frac{\mu_t^{(\gamma^{(\ell)})}(dx)}{dx}.
\end{align}
We have proved the following:

\begin{theorem}{(Stochastic Fokker-Planck equation (II))} \\
Assume that \eqref{abs} holds.
Then $m (t,x)$ satisfies the following SPDE:
\begin{align}
dm(t,x)&= A_0^{*}m(t,x) dt + A_1^{*}m(t,x) dB_1(t);\quad t>0,\nonumber\\
m(0,x) &=\mathcal{L}(X(0)), \label{FP2}
\end{align}
where the operators $A_0^{*},A_1^{*}$ are given by \eqref{A_0m},\eqref{A_1m}, respectively.
\end{theorem}

\subsection {An HJB equation for optimal control of conditional McKean-Vlasov jump diffusions (II): The absolutely continuous case}
In this case the HJB equation obtained in the previous section can be simplified. We outline this in the following:\\
Before we proceed, we introduce and explain some notation:\\
 If $\psi \in C^1(L_K^1(\RR^d))$ then $\nabla_m\psi$ is the Fr\'echet derivative of $\psi$ with respect to $m \in L_K^1(\RR^d)$. Therefore it is a bounded linear functional on $L_K^1(\RR^d)$. By the Riesz representation theorem $\nabla_m\psi$ can be represented by a bounded function $\Psi$ on $\RR^d$, in the sense that the action of $\nabla_m\psi$ on a function $h \in L_K^1(\RR^d)$ can be written
\begin{align*}
\langle \nabla_m \psi, h \rangle = \int_{\RR^d}\Psi(x)h(x) K(dx); \quad h \in L_k^1(\RR^d).
\end{align*}

We consider a controlled version of the system \eqref{MV2}, in which we have introduced a control process $u=\{u(t),t\in[0,T] \}$, i.e. 
\begin{align} \label{MV4}
dX(t)&=dX^{(u)}(t)=\alpha(t,X(t),m(t,\cdot),u(t))dt + \beta(t,X(t),m(t,\cdot),u(t)) dB(t)\\
&+\int_{\mathbb{R}^d} \gamma(t,X(t^-), m(t,\cdot),u(t), \zeta) \widetilde{N}(dt,d\zeta);\quad t>0,\nonumber\\
X(0)&= x \in \mathbb{R}^d,\nonumber
\end{align}
with coefficients $\alpha(t,x,m,u),\beta(t,x,m,u):[0,T]\times \mathbb{R}^d \times L^1(\mathbb{R}^d)\times \mathbb{U} \rightarrow \mathbb{R}$ and $\gamma(t,x,m,u,\zeta):[0,T]\times \mathbb{R}^d\times L_K^1(\mathbb{R}^d)\times \mathbb{U} \times \mathbb{R}^k \rightarrow \mathbb{R}$.\\
Here $m(t,\cdot)$ is the conditional density of $X(t)$ given $\mathcal{F}_t^{(1)}$, in the sense that
\begin{align*}
\int_{\RR^d} g(x)m(t,x)dx=\EE[g(X(t) | \mathcal{F}_t^{(1)}]; \text{ for all } g \text{ such that } \EE[|g(X(t))|] < \infty.
\end{align*}
We say that $u$ is \emph{admissible} if $u$ is Markovian, 
i.e. $u$ has the form $u(t)=u_0(t,X(t),m(t,\cdot))$ for some function $u_0: \RR^3 \mapsto \RR$ (see below), 
and there is a unique solution $X^{(u)}$ satisfying
$\EE[|X^{(u)}(t)|^2] < \infty \text{ for all } t.$
The set of admissible controls is denoted by $\mathcal{A}$. As is customary, for simplicity (and a slight abuse) of notation we do not distinguish in notation between $u_0$ and $u$ in the following.\\
From Theorem 4.3, we have the following Fokker-Planck equation for $m$:
\begin{align}\label{FP3}
dm(t,x)&= A_0^{*,u} m(t,x) dt+ A_1^{*,u} m(t,x) dB_1(t); \quad t>0,\\
m(0,x)&=m_0(x);  \text{ a given initial probability density on } \mathbb{R}^d,\nonumber
\end{align}
where the integro-differential operators $A_0^{*,u},A_1^{*,u}$ associated to the controlled process $u$, are
(see \eqref{A0*},\eqref{A1*}):
\begin{align}
A_0^{*}m&= - \sum_{j=1}^d \frac{\partial}{\partial x_j}[\alpha_j m] +\frac{1}{2}\sum_{n,j=1}^d \frac{\partial^2}{\partial x_j \partial x_n}[(\beta \beta^{(T)})_{n,j} m] \nonumber\\
&+\sum_{\ell=1}^k  \int_{\mathbb{R}}\Big\{m^{(\gamma^{(\ell)})}-m+\sum_{j=1}^d \frac{\partial}{\partial x_j} [\gamma_j^{(\ell)}(s,\cdot,\zeta)m]\Big\} \nu_{\ell} \left( d\zeta \right), \label{A0m}
\end{align}
and
\begin{align}
A_1^{*}m= - \sum_{j=1}^d \frac{\partial}{\partial x_j} [\beta_{1,j} m].\label{A1m}
\end{align}
We introduce now the performance functional:
\begin{align*}
&J_u(s, x,m)=J(y)\nonumber\\
&=\mathbb{E}^{y}\Big[ \int_0^T f(s+t, X(t),m(t,\cdot),u(t))dt  + g(X(T),m(T,\cdot)) \Big],
\end{align*}
where $y=(y_0, y_1,y_2)$. Thus time starts at $y_0=s$, $X(t)$ starts at $x=y_1$ and $m(t)$ starts at $m=y_2$.
Then we define the value function:
\begin{align*}
\Phi(y)=\Phi(s,x,m)= \sup_{u \in \mathcal{A}} J_u (s,x,m),
\end{align*}
where $\mathcal{A}$ denotes the set of admissible controls.\\

We now give sufficient conditions  under which a smooth function $\varphi$ satisfying the HJB equation coincides with the value function $\Phi$, i.e. a verification theorem:

\begin{theorem}\textbf{(An HJB equation for conditional control of McKean-Vlasov jump diffusions (II))} \\
For $v \in \mathbb{U}$ define $G_v$ to be the following integro-differential operator, which is \emph{the generator of $Y(t)=(s+t,X(t),m_t(\cdot))$}  given the control value $v$:
\begin{align*}
G_v \varphi(s,x,m)&= \frac{\partial \varphi}{\partial s} +\sum_{j=1}^d \alpha_j (s,x,m,v)\frac{\partial \varphi}{\partial x_j} + \langle \nabla_{m} \varphi, A_0^{*,v} m \rangle \nonumber\\
&+ \tfrac{1}{2}\sum_{j,n=1}^{d}  (\beta \beta^{T})_{j,n}(s,x,m,v) \frac{\partial ^2 \varphi}{\partial x_j \partial x_n} + \sum_{j=1}^d \beta_{j,1}\frac{\partial}{\partial x_j}\langle\nabla_{m} \varphi,A_1^{*,v}m\rangle \nonumber\\
&+\tfrac{1}{2} 
\langle A_1^{*,v}m, \langle D_{m}^2 \varphi,A_1^{*,v}m\rangle \rangle \nonumber\\
&+\sum_{\ell =1}^k \int_{\mathbb{R}} \{ \varphi(s, x+\gamma^{(\ell)}, m)) - \varphi(s,x,m) 
-\sum_{j=1}^d\gamma_j^{(\ell)}  \tfrac{\partial}{\partial x_j} \varphi(s,x,m) \}\nu_{\ell}(d\zeta) ;\nonumber\\
& \varphi \in C^{1,2,2}([0,T] \times \mathbb{R}^d \times L_K^1(\mathbb{R}^d)),
\end{align*}
where $\gamma=\gamma(s,x,\mu,u,\zeta)$ and $\gamma^{(\ell)}$ is column number $\ell$ of the $d \times k$- matrix $\gamma$.
Suppose there exists a function $\widehat \varphi(s,x,m)\in C^{1,2,2}([0,T] \times \mathbb{R}^d \times L_K^1(\RB^d))$ and a Markov control $\widehat{u}=\widehat{u}(y) \in \AC$ such that, for all  $y$,

\begin{align}
\sup_{v \in \RR} \Big\{f(y,v)+G_v\widehat \varphi(y) \Big \} &= f(y,\widehat{u}(y))+G_{\widehat{u}(y)}\widehat \varphi(y) = 0, \label{HJB2}\\
& \text{ and}\nonumber\\
&\widehat \varphi(T,x,m)=g(x,m).
\end{align}
Then $\widehat{u}$ is an optimal control and $\widehat \varphi=\Phi$.
\end{theorem}

\dproof
The proof is the same as the proof of Theorem 4.1 and is omitted.
\fproof

\section{Examples}
In this section we illustrate our results by solving explicitly some conditional control problems for McKean-Vlasov jump diffusions. For simplicity and without loss of generality, we consider the one dimensional case.
\subsection{A conditional linear-quadratic mean-field control problem} 
 Conditional mean-field linear quadratic control problems have been studied by Pham \& Wei \cite{PW}.
The following example is in some sense a special case of the example in Section 5 of \cite{PW}, except that Pham \& Wei do not consider jumps. Moreover, our approach is different: We use our HJB equation (not square completion) to find an explicit solution.\\

Suppose the system is given by
\begin{align*}
dX\left( t\right) &=u\left( t\right) dt+\theta dB_1(t)\nonumber\\
&+\sigma \mathbb{E}\left[ X( t) |\mathcal{F}_t^{(1)}
\right] dB_2 \left( t\right) +\int_{\mathbb{R}} \gamma_0(\zeta) \EE[X(t)|\mathcal{F}_t^{(1)}] \widetilde{N}(dt,d\zeta);\quad t>0,\\
X(0)&=x \in \mathbb{R},
\end{align*}
with performance functional%
\[
J\left( u\right) =\mathbb{E}\left[ -\frac{1}{2}\int_{0}^{T}u^{2}\left( t\right) dt-%
\frac{1}{2}X^{2}(T) \right]; \quad u \in \mathcal{A},
\]%
where $\theta>0, \sigma >0$ and $\gamma_0(\zeta)$ are a given constants and a given function, respectively. \\
Note that if we define $$q\left( x\right) =x,$$ then $\mathbb{E}\left[ X\left( t\right)|\mathcal{F}_t^{(1)} \right]
=\left\langle \mu_t,q\right\rangle =\int_{\mathbb{R}}x \mu_t(dx)$.\\
In this case we have $\beta=(\theta, \sigma \EE[X(t)|\mathcal{F}_t^{(1)}]) =(\theta, \sigma \langle \mu_t,q\rangle )$ and hence
$$\beta \beta^{(T)}= \theta^2 + \sigma^2 \langle \mu_t,q\rangle ^2.$$
Hence, setting $\gamma_0(\zeta) \EE[X(t)|\mathcal{F}_t^{(1)}]= \gamma$ for simplicity of notation, we see that in this case the operators 
are
\begin{align*}
A_0^{*}\mu=-v D\mu +\frac{1}{2}(\theta^2 + \sigma^2 \langle \mu,q\rangle ^2) D^2 \mu 
+ \int_{\mathbb{R}}\{\mu^{(\gamma)}-\mu+\gamma D \mu \} \nu \left( d\zeta \right),
\end{align*}
and
\begin{align*}
A_1^{*}\mu= - \theta D \mu,
\end{align*}
whose dual operators are, respectively,
\begin{align*}
A_0\mu=v D\mu +\frac{1}{2}(\theta^2 + \sigma^2 \langle \mu,q\rangle ^2) D^2 \mu+ \int_{\mathbb{R}}\{\mu^{(-\gamma)}-\mu-\gamma D \mu \} \nu \left( d\zeta \right),
\end{align*}
and
\begin{align*}
A_1\mu= \theta D \mu,
\end{align*}
where, for notational simplicity, $A_n^{v}=A_n,A_n^{*,v}=A_n^{*}; n=0,1.$\\
Note that $\left\langle \mu_t,q^{\prime }\right\rangle =1,\left\langle
\mu_t,q^{\prime \prime }\right\rangle =0$ with this $q.$\\
Therefore the $HJB$ equation becomes
\begin{align*}
&\sup_{v \in \RR}\Big\{ -\frac{1}{2}v^{2}+\frac{\partial \varphi }{\partial s}%
+v\frac{\partial \varphi }{\partial x}+\left\langle \nabla _{\mu}\varphi ,A_0^{\ast}\mu\right\rangle\nonumber\\
&+\tfrac{1}{2}(\theta^2+\sigma ^{2}\left\langle
\mu,q\right\rangle ^{2})\frac{\partial ^{2}\varphi }{\partial x^{2}}+ \theta \frac{\partial}{\partial x}\langle\nabla_{\mu} \varphi , A_1^{*} \mu\rangle \nonumber\\
&+\tfrac{1}{2} 
\langle A_1^{*}\mu, \langle D_{\mu}^2 \varphi,A_1^{*}\mu\rangle \rangle \nonumber\\
&+\int_{\mathbb{R}} \{ \varphi(s, x+\gamma, \mu)) - \varphi(s,x,\mu) -\gamma  \tfrac{\partial}{\partial x} \varphi(s,x,\mu) \}\nu(d\zeta) \Big\} =0.
\end{align*}
As a candidate for the value function we try a function of the form%
\begin{align}
\varphi \left( s,x,\mu \right) &=\kappa_0(s)+\kappa_{1}\left( s\right) x^{2}+\kappa_{2}\left( s\right)
x\left\langle \mu,q\right\rangle +\kappa_{3}\left( s\right) \left\langle
\mu,q\right\rangle ^{2} \nonumber\\
&=\kappa_0(s)+ \kappa_{1}\left( s\right) x^{2}+\kappa_{2}\left( s\right)
x \psi_1(\mu) +\kappa_{3}\left( s\right) \psi_2(\mu) , \label{phi}
\end{align}
where $\kappa_0, \kappa_{1},\kappa_{2},\kappa_{3}$ are deterministic, differentiable functions and we have defined
\begin{align*}
\psi_1(\mu)= \langle \mu,q \rangle \quad \text{ and } \quad \psi_2(\mu)=\langle \mu,q \rangle ^2; \quad \mu \in \mathbb{M}.
\end{align*}
Then, since $\psi_1$ is linear (see Appendix), 
\begin{align*}
\langle \nabla_{\mu} \psi_1, A_0^{*}\mu \rangle  = \langle \psi_1, A_0^{*} \mu \rangle =\langle A_0^{*}\mu,q \rangle = \langle \mu, A_0 q \rangle = v \langle \mu,q^{\prime} \rangle = v,
\end{align*}
and, by the chain rule (see Appendix),
\begin{align*}
\langle \nabla_{\mu} \psi_2, A_0^{*}\mu \rangle = 2 \langle \mu,q \rangle \langle \nabla_{\mu}\psi_1, A_0^{*}\mu \rangle= 2 \langle \mu,q \rangle \langle \mu, A _0 q \rangle &= 2v \langle \mu,q \rangle.
\end{align*}
Therefore 
\begin{equation*}
\left\langle \nabla _{\mu}\varphi ,A_0^{\ast }\mu\right\rangle =v [\kappa_{2}(s)x +2\kappa_{3}(s) \left\langle \mu,q\right\rangle].
\end{equation*}
Similarly,
\begin{align*}
\langle \nabla_\mu \psi_1,A_1^{*}\mu\rangle&= \langle \psi_1,A_1^{*}\mu\rangle= \langle A_1^{*}\mu,q\rangle=\langle\mu,A_1 q\rangle =\langle \mu, \theta q^{\prime}\rangle=\theta,\\
\langle \nabla_\mu \psi_2, A_1^{*}\mu\rangle&= 2 \langle \mu,q \rangle \langle \nabla_\mu \psi_1, A_1^{*}\mu \rangle= 2 \langle \mu,q \rangle \langle \mu, A _1 q \rangle = 2\theta \langle \mu,q \rangle,\nonumber\\
\frac{\partial }{\partial x} \langle \nabla_\mu \varphi, A_1^{*}\mu\rangle&= \frac{\partial}{\partial x}[\theta \kappa_2(s) x] =\theta\kappa_2(s). 
\end{align*}
Since $D_{\mu}^2 \psi_1=0$ and $D_{\mu}^2 \psi_2 (h,k)= 2hk$ (see Appendix), we get
 $D_{\mu}^2\varphi (h,k)= 2\kappa_3(s)hk$, and hence
\begin{align*}
\langle A_1^{*}\mu, \langle D_{\mu}^2 \varphi, A_1^{*}\mu\rangle\rangle= 2 \kappa_3(s)\langle A_1^{*}\mu,q\rangle\langle A_1^{*}\mu,q\rangle=2\theta^2 \kappa_3(s).
\end{align*}
Hence the HJB \eqref{HJB1} becomes, with $\kappa_{i}^{\prime }=\kappa_{i}^{\prime
}\left( s\right) =\frac{d}{ds}\kappa_{i}\left( s\right),$%
\begin{align}
&\sup\limits_{v}\Big\{ -\frac{1}{2}v^{2}+[\kappa_0^{\prime}+ \kappa_{1}^{\prime
}x^{2}+\kappa_{2}^{\prime }x\langle \mu,q\rangle +\kappa_{3}^{\prime
}\langle \mu,q\rangle ^{2}] +v [ 2\kappa_{1}x+\kappa_{2}\langle \mu,q\rangle]  \nonumber \\
&+ v [ x\kappa_{2}+2\kappa_{3}\langle \mu,q\rangle ]+\frac{1}{2}%
(\theta^2+ \sigma ^{2}\langle \mu,q\rangle ^{2})2\kappa_{1} \nonumber\\
&+\theta^2 (\kappa_2+\kappa_3) +\kappa_1 \int_{\mathbb{R}} \gamma^2(\zeta) \nu(d\zeta)
\Big\}=0 .  \nonumber
\end{align}
Maximising with respect to $v$ gives the first order condition%
\[
-v+\left[ 2\kappa_{1}x+\kappa_{2}\left\langle \m,q\right\rangle
+ \kappa_{2}x+ 2\kappa_{3}\left\langle \m,q\right\rangle \right] =0,
\]%
or%
\begin{equation}\label{(3)}
v=\widehat{u}=\left( 2\kappa_{1}+\kappa_{2}\right) x+\left( \kappa_{2}+2\kappa_{3}\right)
\left\langle \m,q\right\rangle.   
\end{equation}%
It remains to verify that with this value of $u=\widehat{u}$, we get \eqref{(1)} (without the $sup$) satisfied:

Substituting \eqref{(3)} into \eqref{(1)} we get%
\begin{align}
&-\frac{1}{2}\Big[ \left( 2\kappa_{1}+\kappa_{2}\right) ^{2}x^{2}+2\left(
2\kappa_{1}+\kappa_{2}\right) \left( \kappa_{2}+2\kappa_{3}\right) x\left\langle
\mu,q\right\rangle +\left( \kappa_{2}+2\kappa_{3}\right) ^{2}\left\langle
\mu,q\right\rangle ^{2}\Big]   \nonumber \\
&+\Big[ \kappa_0^{\prime} +\kappa_{1}^{\prime }x^{2}+\kappa_{2}^{\prime }x\left\langle \mu,q\right\rangle
+\kappa_{3}^{\prime }\left\langle \mu,q\right\rangle ^{2}\Big]   \nonumber \\
&+\Big[ \left( 2\kappa_{1}+\kappa_{2}\right) x+\left( \kappa_{2}+2\kappa_{3}\right)
\left\langle \mu,q\right\rangle \Big] \Big[ 2\kappa_{1}x+\kappa_{2}\left\langle
\mu,q\right\rangle \Big] \nonumber\\
&+(\theta^2+\sigma ^{2}\left\langle
\mu,q\right\rangle ^{2})\kappa_{1} +\theta^2 (\kappa_2+\kappa_3) +\kappa_1 \int_{\mathbb{R}} \gamma^2  \nu(d\zeta) \nonumber \\
&+\Big[ \left( 2\kappa_{1}+\kappa_{2}\right) x+\left( \kappa_{2}+2\kappa_{3}\right)
\left\langle \mu,q\right\rangle \Big] \Big[ \kappa_{2}x+2\kappa_{3}\left\langle
\mu,q\right\rangle \Big]   \nonumber \\
&=\kappa_0^{\prime}+ x^{2}\Big[ -\frac{1}{2}\left( 2\kappa_{1}+\kappa_{2}\right) ^{2}+\kappa_{1}^{\prime
}+\left( 2\kappa_{1}+\kappa_{2}\right) 2\kappa_{1}+\left( 2\kappa_{1}+\kappa_{2}\right) \kappa_{2}\Big] \nonumber\\
&+x\left\langle \mu,q\right\rangle \Big[ -\left( 2\kappa_{1}+\kappa_{2}\right) \left(
\kappa_{2}+2\kappa_{3}\right) +\kappa_{2}^{\prime }+\left( \kappa_{2}+2\kappa_{3}\right)
2\kappa_{1}+\left( \kappa_{2}+2\kappa_{3}\right) \kappa_{2}\Big]   \nonumber \\
&+\left\langle \mu,q\right\rangle^2 \Big[ -\frac{1}{2}\left(
\kappa_{2}+2\kappa_{3}\right) ^{2}+\kappa_{3}^{\prime }+\kappa_{2}\left( \kappa_{2}+2\kappa_{3}\right)\nonumber\\
&+\kappa_{1}\sigma ^{2}+\int_{\mathbb{R}} \gamma_0^2 \nu(d\zeta))+2\kappa_{3}\left( \kappa_{2}+\kappa_{3}\right) \Big] +\theta^2(\kappa_1 + \kappa_2+\kappa_3) =0.
\nonumber
\end{align}
This holds for all $x, \left\langle \mu,q\right\rangle $ if and only if the
following system of differential equations (Riccati equations) are satisfied%
\begin{align}
&\kappa_0^{\prime}+ \theta^2(\kappa_1+\kappa_2+\kappa_3)=0;\quad \kappa_0(T)=0, \label{R0}\\
&\kappa_{1}^{\prime }+\frac{1}{2}\left( 2\kappa_{1}+\kappa_{2}\right) ^{2}=0;\quad \kappa_{1}\left( T\right) =-%
\frac{1}{2},  \label{R1}
\\
&\kappa_{2}^{\prime }+2\left( 2\kappa_{1}+\kappa_{2}\right) \left( \kappa_{2}+2\kappa_{3}\right)
=0
;\quad \kappa_{2}\left( T\right) =0,  \label{R2}\\
&\kappa_{3}^{\prime }+\frac{1}{2}\left( \kappa_{2}+2\kappa_{3}\right) ^{2} + \kappa_1 (\sigma
^{2}+\int_{\mathbb{R}} \gamma^2 \nu(d\zeta))=0;\quad \kappa_{3}\left( T\right) =0.  \label{R3}
\end{align}
By general theory of matrix Riccati equations (see e.g. Theorem 37 in \cite{S}) there is a unique solution 
$\widehat{\kappa}_0(s), \widehat{\kappa}_{1}\left( s\right) ,\widehat{\kappa}_{2}\left( s\right) ,%
\widehat{\kappa}_{3}\left( s\right) $ satisfying (\ref{R0}), (\ref{R1}), (\ref{R2}), (\ref%
{R3}), respectively. \\
We conclude that, with these choices of $\kappa_j=\widehat{\kappa}_j; j=0,1,2,3,$ our guessed candidate $\varphi$ given by \eqref{phi} satisfies all the conditions of the HJB equations, and hence it is indeed the value function. Therefore we have the following result:

\begin{theorem}
The optimal control $\ \widehat{u}\left( t\right) $ of the problem to
maximize%
\[
J\left( u\right) :=\EE\left[ -\frac{1}{2}\int_{0}^{T}u^{2}\left( t\right) dt-%
\frac{1}{2}X^{2}\left( T\right) \right]; \quad u \in \mathcal{A},
\]%
is given in feedback form as follows:%
\begin{equation*}
\widehat{u}\left( Y(t)\right) =\Big( 2\widehat{\kappa}_{1}\left( t\right) +%
\widehat{\kappa}_{2}\left( t\right) \Big) X\left( t\right) +\Big( \widehat{\kappa}%
_{2}\left( t\right) +2\widehat{\kappa}_{3}\left( t\right) \Big) \EE\left[ X\left(
t\right)|\mathcal{F}_t^{(1)} \right] . 
\end{equation*}%
The value function is%
\begin{align*}
&\Phi(s,x,\mu)=\widehat{\varphi }\left( s,x,\mu\right) :=\widehat{\kappa}_0(s)+\widehat{\kappa}_{1}\left( s\right) x^{2}+%
\widehat{\kappa}_{2}\left( s\right) x\left\langle \mu,q\right\rangle +\widehat{\kappa}%
_{3}\left( s\right) \left\langle \mu,q\right\rangle ^{2}\nonumber\\
&=\widehat{\kappa}_0(s)+\widehat{\kappa}_{1}\left( s\right) x^{2}+%
\widehat{\kappa}_{2}\left( s\right) x\int_{\RR} z\mu(dz)  +\widehat{\kappa}%
_{3}\left( s\right) (\int_{\RR} z \mu(dz))^2\nonumber\\
&=\widehat{\kappa}_0(s)+\widehat{\kappa}_{1}\left( s\right) x^{2}+%
\widehat{\kappa}_{2}\left( s\right) x \EE[X(s)|\mathcal{F}_s^{(1)} ]+\widehat{\kappa}%
_{3}\left( s\right) (\EE[X(s)|\mathcal{F}_s^{(1)} ] )^2.
\end{align*}
\end{theorem}
\subsection{An optimal consumption/harvesting problem}
The following problem may for example be considered as an optimal fish population harvesting problem, or as a problem of optimal consumption from a cash flow.

Consider the following controlled conditional mean-field SDE:%
\begin{align*}
dX\left( t\right) &= \left( \rho \left( t\right) -c\left( t\right) \right) \EE
\left[ X\left( t\right)|\mathcal{F}_t^{(1)} \right] dt\nonumber\\
&+\theta \EE[X(t)|\mathcal{F}_t^{(1)}]dB_1(t)+\sigma _{0}\left( t\right) \EE\left[
X\left( t\right)|\mathcal{F}_t^{(1)} \right] dB_2\left( t\right) \nonumber\\
&+\int_{\mathbb{R}}\gamma
_{0}\left( t,\zeta \right) \EE\left[ X\left( t\right) |\mathcal{F}_t^{(1)}\right] \widetilde{N}%
\left( dt,d\zeta \right) ;\quad t>0, \nonumber\\ 
X\left( 0\right) = & x,
\end{align*}

where $\rho \left( t\right) ,\sigma_0 \left( t\right) $ and $\gamma_0(t,\zeta)$ are bounded
deterministic functions. The consumption rate $c =c\left( t,\omega\right); (t,\omega)\in \ [0,T] \times \Omega $ is defined to be an $\mathbb{F}-$predictable positive process; it is admissible if $\mathbb{E}[\int_0^Tc(t)^2dt]<\infty$.

We want to maximize the total expected utility form the consumption,
expressed by the performance functional%
\begin{equation}
J\left( c\right) =\EE\left[ \int_{0}^{T}\ln \left( c\left( t\right) \EE\left[
X\left( t\right)|\mathcal{F}_t^{(1)} \right] \right) dt+\lambda \ln\Big(\EE\left[ X\left(T\right)|\mathcal{F}_T^{(1)} \right] \Big)
\right];\quad c \in \mathcal{A} ,  \label{(2)}
\end{equation}%
where $\lambda >0$ is a (deterministic) constant.

To put this into our framework, we proceed as in the previous example and write, with $\mu_t=\mu(t, \cdot)$,
\[
\EE\left[ X\left( t\right)|\mathcal{F}_t^{(1)} \right] =\left\langle \m_{t},q\right\rangle \text{
where }q\left( x\right) =x.
\]%
Define the operator $A_0$ by
\begin{align*}
A_0\m(x)&=\left( \rho \left( t\right) -c\left( t\right) \right) \left\langle
\m,q\right\rangle D\m+\tfrac{1}{2}(\theta^2+\sigma_0 ^{2}\left( t\right)
\left\langle \m,q\right\rangle ^{2}) D^2 \m \nonumber\\
&+ \int_{\mathbb{R}}\{\mu^{(-\gamma)}-\mu-\gamma D \mu \} \nu \left( d\zeta \right) ;\quad \m\in \mathbb{M}, 
\end{align*}
where $\gamma= \gamma_0 \langle\m,q\rangle$.
The adjoint of this operator is
\begin{align}
A_0^{*}\m&=-\left( \rho \left( t\right) -c\left( t\right) \right) \left\langle
\m_{t},q\right\rangle D\m+\tfrac{1}{2}(\theta^2+\sigma_0 ^{2}\left( t\right)
\left\langle \m_{t},q\right\rangle ^{2}) D^2 \m \nonumber\\
&+ \int_{\mathbb{R}}\{\mu^{(\gamma)}-\mu+\gamma D \mu \} \nu \left( d\zeta \right) ;\quad \m\in \mathbb{M}.  \label{(4)}
\end{align}
Similarly, we define
\begin{align*}
A_1 \m&= \theta \l \m,q \r D\m,  \text{  whose adjoint is  } \\
A_1^{*}\mu&= -\theta \<\m,q\>  D\m; \quad \m \in \mathbb{M}.
\end{align*}
The corresponding HJB equation for the value function $\varphi $ becomes%
\begin{align}\label{(1)}
&\sup\limits_{c > 0}\Big\{ \ln \left( c \langle
\m,q\rangle \right) +\frac{\partial \varphi }{\partial t}+\left( \rho
\left( t\right) -c\right) \langle \m,q \rangle \frac{\partial
\varphi }{\partial x}+\langle \nabla _{\m}\varphi ,A_0^{\ast }\m \rangle\nonumber\\
&+\tfrac{1}{2}(\theta^2+\sigma_0 ^{2}\left( t\right) \langle
\m_{t},q \rangle ^{2})\frac{\partial ^{2}\varphi }{\partial x^{2}}  +\theta \frac{\partial}{\partial x} \<\nabla_{\m} \varphi, A_1 ^{*}\m\> \nonumber\\
&+\tfrac{1}{2} \theta^2 \<\m,q\>^2
\langle A_1^{*}\mu, \langle D_{\mu}^2 \varphi,A_1^{*}\mu\rangle \rangle\nonumber\\
&+\int_{\mathbb{R}} \{ \varphi(s, x+\gamma, \m)) - \varphi(s,x,\m) -\gamma  \tfrac{\partial}{\partial x} \varphi(s,x,\m) \}\nu(d\zeta) \Big\}=0
;\quad t<T,
\end{align}
with terminal value%
\begin{equation}
\varphi \left( T,x,\m\right) =\lambda \ln \left\langle \m_{T},q\right\rangle. 
\label{(6)}
\end{equation}%
Let us guess that the value function is of the form%
\begin{equation}
\varphi \left( s,x,\m\right) =\kappa_{0}\left( s\right) +\kappa_{1}\left( s\right) \ln
\left\langle \m,q\right\rangle ,  \label{(7)}
\end{equation}%
for some $C^{1}$ deterministic functions $\kappa_{0}\left( s\right) ,\kappa_{1}\left(
s\right) .$
Define
\begin{align*}
\psi(\m)=\ln \langle \m, q \rangle.
\end{align*}
Then we have (see Appendix), 
\begin{align*}
\nabla _{\mu}\psi(\m) (h) &= \frac{\<h,q\>}{\<\m,q\>}\\
D^2\psi(\m)(h,k)&=-\frac{\<h,q\> \<k,q\>}{\<\m,q\>^2}.
\end{align*}
Hence
\small
\begin{align*}
\langle \nabla _{\m}\varphi ,A_0^{\ast }\m\rangle  &=\kappa_{1}\left(
s\right) \frac{1}{\left\langle \m,q\right\rangle }\left\langle \nabla
_{\m}\psi ,A_0^{\ast }\m\right\rangle \nonumber\\
&=\kappa_{1}\left( s\right) \frac{1}{\left\langle \m,q\right\rangle }\left\langle
A_0^{\ast }\m,q\right\rangle =\kappa_{1}\left( s\right) \frac{1}{\left\langle
\m,q\right\rangle }\left\langle \m,A_0 q\right\rangle   \nonumber \\
&=\kappa_{1}\left( s\right) \left( \rho -c\right) \left[ \left\langle
\m,q^{\prime }\right\rangle +\left\langle \m,\int_{\mathbb{R}}\left\{ q\left( x-\gamma_0\langle \m,q\rangle
\right) -q+\gamma_0 \langle \m,q\rangle q^{\prime }\right\} \nu \left( d\zeta \right)
\right\rangle \right]   \nonumber \\
&=\kappa_{1}\left( s\right) \left( \rho \left( s\right) -c\right).  
\end{align*}%
Here we have used that

\begin{align*}
& \langle \m,q^{\prime }\rangle =\left\langle \m,1\right\rangle =\int_{\mathbb{R}}\m\left(x\right) dx=1 \quad \text{ and, since } q \text{ is linear }, \\ 
&\int_{\mathbb{R}} \{\left\langle \m,q\left( x-\gamma \right) \right\rangle
-\left\langle \m,q\right\rangle +\gamma \left\langle \m,q^{\prime
}\right\rangle\} \nu \left( d\zeta \right)   
=\int_{\mathbb{R}}\left\{ \left\langle \m,q\right\rangle -\left\langle
\m,\gamma \right\rangle -\left\langle \m,q\right\rangle +\gamma \right\} \nu
\left( d\zeta \right) =0.%
\end{align*}
Similarly, we get
\begin{align*}
\frac{\partial}{\partial x} \langle\nabla_\m \varphi, A_1^{*} \m\rangle =
0
\end{align*} 
and
\begin{align*}
\langle A_1^{*}\mu, \langle D_{\mu}^2 \varphi,A_1^{*}\mu\rangle \rangle
&=-\frac{\kappa_1(s)\<A_1^{*}\m,q\> \<A_1^{*}\m,q\>}{\<\m,q\>^2}\nonumber\\
&=-\frac{\kappa_1(s)\<\m,A_1 q\> \<\m,A_1 q\>}{\<\m,q\>^2}=-\frac{\kappa_1(s)\theta^2}{\<\m,q\>^2}.
\end{align*}
Therefore the HJB equation now takes the form%
\begin{align}
&\sup\limits_{c > 0}\Big\{ \ln c+\ln \left\langle
\m,q\right\rangle +\kappa_{0}^{\prime }\left( s\right)+\kappa_{1}^{\prime }\left(s\right) \ln \langle \m,q \rangle
+\kappa_{1}\left( s\right) \left( \rho \left( t\right) -c\right) -\tfrac{1}{2} \kappa_1(s) \theta^4\Big\}
=0.
\end{align}
\label{(9)}
The optimising value of $c$ is the solution of the equation%
\[
\frac{1}{c}-\kappa_{1}\left( s\right) =0,
\]%
i.e%
\begin{equation*}
c\left( s\right) =\widehat{c}\left( s\right) =\frac{1}{\kappa_{1}\left( s\right) }.
\end{equation*}%
Substituting this into \eqref{(9)} we get%
\begin{equation}
\begin{array}{l}
-\ln \kappa_{1}\left( s\right) +\ln \left\langle \m,q\right\rangle
+\kappa_{0}^{\prime }\left( s\right) +\kappa_{1}^{\prime }\left( s\right) \ln
\left\langle \m,q\right\rangle  \\ 
+\rho \left( s\right) \kappa_{1}\left( s\right)+\kappa_1(s) \theta -\tfrac{1}{2}\kappa_1(s)\theta^4=0.%
\end{array}
\label{(11)}
\end{equation}%
From (\ref{(6)}) and (\ref{(7)}) we get the terminal values%
\begin{equation}
\kappa_{0}\left( T\right) =0,\quad \kappa_{1}\left( T\right) =\lambda.   \label{(12)}
\end{equation}%
Hence, if we choose $\kappa_{1}$ such that%
\begin{eqnarray}
\kappa_{1}^{\prime }\left( s\right)  &=&-1,\text{ i.e}  \label{(13)} \\
\kappa_{1}\left( s\right)  &=&\lambda +T-s,  \nonumber
\end{eqnarray}%
and let $\kappa_{0}\left( s\right) $ be the given by%
\begin{equation}
\left\{ 
\begin{array}{ll}
\kappa_{0}^{\prime }\left( s\right) = & \tfrac{1}{2}\kappa_1(s)  \theta^4 +\ln \kappa_{1}\left( s\right) -(\rho \left(
s\right) +\theta)\kappa_{1}\left( s\right); s<T, \\ 
\kappa_{0}\left( T\right) = & 0,
\end{array}%
\right.   \label{(14)}
\end{equation}%
we see that (\ref{(11)}) holds. 
We have proved the following:
\begin{theorem}

The function 
\[
\varphi \left( s,x,\m\right) =\kappa_{0}\left( s\right) +\kappa_{1}\left( s\right) \ln
\left\langle \m,q\right\rangle,
\]%
with $\kappa_{0},\kappa_{1}$ defined by (\ref{(13)}), (\ref{(14)}), satisfies all the
conditions of our HJB theorem and therefore
\[
\Phi \left( s,x,\m\right) =\varphi \left( s,x,\m\right) 
\]%
is the value function, and%
\[
\widehat{c}\left( s\right) =\frac{1}{\kappa_{1}\left( s\right) }
\]%
is the optimal control.
\end{theorem}

\section{Summary}
\begin{itemize}
 \item In this paper we study the  \emph{conditional McKean-Vlasov jump diffusions}, which are conditional mean-field stochastic differential equations with jumps.
 \item Using Fourier transforms of Radon measures we prove that the conditional law process of the solution of a McKean-Vlasov jump diffusion satisfies (in the sense of distributions)  a stochastic Fokker-Planck equation.
  \item Combining the Fokker-Planck equation with the original equation we represent the solution of the conditional McKean-Vlasov jump diffusion as the solution of a multidimensional  Markovian system. This allows us to formulate an HJB equation for the optimal control of such systems.
\item Finally we illustrate our results by solving explicitly some optimal control problems for conditional McKean-Vlasov jump diffusions.
 \end{itemize}

\section{Appendix: Double Fr\' echet derivatives}
In this section we recall some basic facts we are using about the Fr\' echet derivatives of a function $f: V \mapsto W$, where $V,W$ are given Banach spaces.
\begin{definition}
We say that $f$ has a Fr\' echet derivative $\nabla_xf=Df(x)$ at $x \in V$ if there exists a bounded linear map $A:V \mapsto W$ such that
\begin{align*}
\lim_{h \rightarrow 0} \frac{\|f(x+h)-f(x)-A(h)\|_{W}}{\|h\|_{V}} =0
\end{align*}
Then we call $A$ the Fr\' echet derivative of $f$ at x and we put $Df(x) =A$
\end{definition}
Note that $Df(x) \in L(V,W)$ (the space of bounded linear functions from $V$ to $W$),  for each $x$.

\begin{definition}
We say that $f$ has a double Fr\' echet derivative $D^2 f(x)$ at $x$ if there exists a bounded bilinear map $A(h,k): V \times V \mapsto W$ such that
\begin{align*}
\lim_{k \rightarrow 0} \frac{\|Df(x+k)(h)-Df(x)(h)-A(h,k)\|_{W}}{\|h\|_{V}} =0
\end{align*}
\end{definition}

\begin{example}
\begin{itemize}
\item
Suppose $f:\mathbb{M} \mapsto \mathbb{R}$ is given by
\begin{align*}
f(\mu)=\<\mu,q\>^2, \text{ where } q(x)=x.
\end{align*}
Then 
\begin{align*}
f(\mu +h) -f(\mu)&= \<\mu+h,q\>^2 -\<\mu,q\>^2\nonumber\\
&= 2 \<\mu,q\> \<h,q\>+\<h,q\>^2, 
\end{align*}
so we see that
\begin{equation*}
Df(\mu)(h)=2\<\mu,q\>\<h,q\>.
\end{equation*}
To find the double derivative we consider
\begin{align*}
&Df(\mu+k)(h)-Df(\mu)(h)\nonumber\\
&=2\<\mu+k,q\>\<h,q\>-2\<\mu,q\>\<h,q\>\nonumber\\
&=2\<k,q\>\<h,q\>,
\end{align*}
and we conclude that
\begin{equation*}
D^2f(\mu)(h,k)=2\<k,q\>\<h,q\>.
\end{equation*}
\item
Next, assume that $g:\mathbb{M} \mapsto \RR$ is given by $g(\mu)=\<\mu,q\>$.
Then, proceeding as above we find that
\begin{align*}
Dg(\mu)(h)&=\<h,q\> \text{ (independent of } \mu)\\
&\text{ and }\nonumber\\
D^2g(\mu) &=0.
\end{align*}

\end{itemize}
\end{example}
\thanks{{\bf Acknowledgements}. Many thanks to the editor and the anonymous referees for their insightful and valuable comments that helped improve the content of the paper.}



\begin{thebibliography}{99}

\bibitem{A} Agram, N. (2019). Stochastic optimal control of McKean-Vlasov equations with anticipating law. Afrika Matematika, 30(5), 879-901.

\bibitem{AHO}
Agram, N., Hu, Y., \& \O ksendal, B. (2022). Mean-field backward stochastic differential equations and applications. Systems \& Control Letters, 162, 105196.

\bibitem{AD} Andersson, D. \& Djehiche, B. (2011). A maximum principle for SDEs of mean-field type. Applied Mathematics \& Optimization, 63(3), 341-356.


\bibitem{BR} Barbu, V., \& R\" ockner, M. (2021). Solutions for nonlinear Fokker-Planck equations with measures as initial data and McKean-Vlasov equations. Journal of Functional Analysis, 280(7), 108926.

\bibitem{BR1} Barbu, V., \& R\" ockner, M. (2021). Uniqueness for nonlinear Fokker-Planck equations and weak uniqueness for McKean-Vlasov SDEs. Stochastics and Partial Differential Equations: Analysis and Computations, 9(3), 702-713.

\bibitem{BC} Bayraktar, E., \& Chakraborty, P. (2021). Mean field control and finite dimensional approximation for regime-switching jump diffusions. arXiv preprint arXiv:2109.09134.

\bibitem{BHL} Bensoussan, A., Huang, T., \& Lauri\`ere, M. (2018). Mean field control and mean field game models with several populations. Minimax Theory and its Applications, 03 (2), 173--209.

\bibitem{BKRS} Bogachev, V. I., Krylov, N. V., R\" ockner, M., \& Shaposhnikov, S. V. (2015). Fokker-Planck-Kolmogorov Equations (Vol. 207). American Mathematical Society.

\bibitem{BLPR} Buckdahn, R., Li, J., Peng, S., \& Rainer, C. (2017). Mean-field stochastic differential equations and associated PDEs. The Annals of Probability, 45(2), 824-878.

\bibitem{BIRS} Burzoni, M., Ignazio, V., Reppen, A. M., \& Soner, H. M. (2020). Viscosity solutions for controlled McKean--Vlasov jump-diffusions. SIAM Journal on Control and Optimization, 58(3), 1676-1699.

\bibitem{C} Cardaliaguet, P. (2010). Notes on mean field games (p. 120). Technical report.

\bibitem{CD} Carmona, R., Delarue, F., \& Lachapelle, A. (2013). Control of McKean-Vlasov dynamics versus mean field games. Mathematics and Financial Economics, 7(2), 131-166.

\bibitem{CG} Coghi, M., \& Gess, B. (2019). Stochastic nonlinear fokker-planck equations. Nonlinear Analysis, 187, 259-278.

\bibitem{CGKPR} Cosso, A., Gozzi, F., Kharroubi, I., Pham, H., \& Rosestolato, M. (2021). Master Bellman equation in the Wasserstein space: Uniqueness of viscosity solutions. arXiv preprint arXiv:2107.10535v2.

\bibitem{F} {Folland, G.B.} (1984). Real Analysis. Modern Techniques and Their Applications. Wiley.

\bibitem{GPW} Guo, X., Pham, H., \& Wei, X. (2020). It\^ o's formula for flow of measures on semimartingales. arXiv preprint arXiv:2010.05288.

\bibitem{H} Hafayed, M. (2013). A mean-field maximum principle for optimal control of forward-backward stochastic differential equations with Poisson jump processes. International Journal of Dynamics and Control, 1(4), 300-315.

\bibitem{JMW} Jourdain, B., M\' el\' eard, S., \& Woyczynski, W. (2008). Nonlinear SDEs driven by L\' evy processes and related PDEs. Alea, 4, 1-29.

\bibitem{KP} Kharroubi, I., \& Pham, H. (2015). Feynman-Kac representation for Hamilton-Jacobi-Bellman IPDE. The Annals of Probability, 43(4), 1823-1865.

\bibitem{KT} Kolokoltsov, V., \& Troeva, M. (2015). On the mean field games with common noise and the McKean-Vlasov SPDEs. arXiv preprint arXiv:1506.04594.

\bibitem{KOR} Kromer, E., Overbeck, L., \& R\" oder, J. A. L. (2015). Feynman-Kac for functional jump diffusions with an application to Credit Value Adjustment. Statistics \& Probability Letters, 105, 120-129.

\bibitem{KX} Kurtz, T. G., \& Xiong, J. (1999). Particle representations for a class of nonlinear SPDEs. Stochastic Processes and their Applications, 83(1), 103-126.

\bibitem{LL} Lasry, J-M. \& Lions, P.-L. (2007) Mean-field games. Japan J. Math. 2, 229-260.

\bibitem{LP} Lauri\`ere, M. \& Pirroneau, O. (2014). Dynamic programming for mean-field type control. Comptes Rendus Mathematique, 352(9), 707-713.

\bibitem{LP1} Lauri\`ere, M. \& Pirroneau, O. (2016). Dynamic programming for mean-field type control. Journal of Optimization Theory and Applications, 169(3), 902-924.


\bibitem{M} McKean Jr, H. P. (1966). A class of Markov processes associated with nonlinear parabolic equations. Proceedings of the National Academy of Sciences of the United States of America, 56(6), 1907.


\bibitem{MP} Miller, E., \& Pham, H. (2019). Linear-quadratic McKean-Vlasov stochastic differential games. In Modeling, Stochastic Control, Optimization, and Applications (pp. 451-481). Springer.

\bibitem{O} \O ksendal, B. (2013). Stochastic Differential Equations. 6th edition. Springer.

\bibitem{OS} \O ksendal, B. \& Sulem, A. (2019). Applied Stochastic Control of Jump diffusions. 3rd edition. Springer.

\bibitem{PW} Pham, H., \& Wei, X. (2017). Dynamic programming for optimal control of stochastic McKean--Vlasov dynamics. SIAM Journal on Control and Optimization, 55(2), 1069-1101.

\bibitem{S} Sontag, E. (1998).  Mathematical Control Theory: Deterministic Finite Dimensional Systems, 2nd Edition Texts in Applied Mathematics, Volume 6, Second Edition, New York: Springer. 



\end{thebibliography}
\end{document}